\documentclass[american]{amsart}

\usepackage[T1]{fontenc}

\usepackage[latin1]{inputenc}

\usepackage{babel}

\newfont{\gothic}{eufm10 scaled 1100}

\newfont{\smgothic}{eufm10 scaled 900}

\makeatletter

\theoremstyle{plain}    

\newtheorem{thm}{Theorem}[section]

\numberwithin{equation}{subsection} 

\numberwithin{figure}{section} 

\theoremstyle{plain}    

\newtheorem{lem}[thm]{Lemma} 

\theoremstyle{plain}    

\newtheorem{keylem}[thm]{Key Lemma} 

\theoremstyle{plain}

\theoremstyle{plain}    

\newtheorem{Ext-thm}[thm]{Extension Theorem} 

\theoremstyle{plain}    

\newtheorem{Def}[thm]{Definition} 

\theoremstyle{plain}    

\newtheorem{prop}[thm]{Proposition} 

\theoremstyle{plain}    

\theoremstyle{plain}    

\theoremstyle{remark}

\theoremstyle{remark}

\newtheorem{notation-assumptions}[thm]{Notation-Assumptions}

\usepackage[all]{xy}

\xyoption{dvips}

\CompileMatrices

\makeatother

\begin{document}

\title{Tsuji's Numerical Trivial Fibrations}


\author{Thomas Eckl}




\address{Thomas Eckl, Institut für Mathematik, Universität
Bayreuth, 95440 Bayreuth, Germany}

\email{thomas.eckl@uni-bayreuth.de}

\urladdr{http://btm8x5.mat.uni-bayreuth.de/\~{}eckl}

\begin{abstract}
The Reduction Map Theorem
in H. Tsuji's work \cite{Tsu00} on 
numerical trivial fibrations is corrected and proven. To this purpose various 
definitions of Tsuji's new intersection 
numbers for pseudo-effective line bundles equipped with a positive singular 
hermitian metric are compared and their equivalence
on sufficiently general smooth curves is shown. Numerically trivial varieties 
are characterized by a decomposition property of the curvature current. An 
important 
adjustment to the Reduction Map Theorem is to 
consider the fact that plurisubharmonic functions are singular on pluripolar 
sets.   
\end{abstract}

\maketitle

\tableofcontents

\bibliographystyle{alpha}


\section{Introduction}

\noindent
In \cite{Tsu00}, H. Tsuji stated several very interesting assertions on the 
structure of pseudo-effective line bundles 
$L$
on a projective manifold
$X$.
In particular he postulated the existence of a meromorphic ``reduction map'',
which essentially says that through the general point of
$X$
there is a maximal irreducible
$L-\!\!$
flat subvariety. Moreover the reduction map should be almost holomorphic, i.e.
has compact fibres which do not meet the indeterminacy locus of the reduction
map. The proofs of \cite{Tsu00} use deep analytic methods and are extremely 
difficult to follow. In \cite{BCEKPRSW00}, the existence of a similar 
reduction map for nef line bundles and that it is almost holomorphic is proven
by purely algebraic methods. This was also stated explicitly in \cite{Tsu00},
but it is not clear how the two reduction maps are connected.

\noindent
The purpose of this note is first to clarify the exact meaning of the 
intersection numbers to which the term ``numerical trivial'' refers. Three
definitions may be found useful in Tsuji's arguments, and it is a very subtle
question when they are equivalent:
\begin{Def} \label{Int-def}
Let
$X$
be a smooth projective complex manifold, let
$L$
be a holomorphic line bundle on
$X$
with positive singular hermitian metric
$h$.
If
$C \subset X$
is an irreducible curve with normalization
$\pi: \tilde{C} \rightarrow C$
such that
$h$ 
is well defined on
$C$,
i.e. 
$h_{|C} \not\equiv + \infty$,
then define the intersection number
\[ (L,h).C := \limsup_{m \rightarrow \infty} \frac{1}{m} 
              h^0(\tilde{C}, \mathcal{O}_{\tilde{C}}(m\pi^\ast L) \otimes 
                             \mathcal{I}((\pi^\ast h)^m)). \]
Here,
$\mathcal{I}((\pi^\ast h)^m))$
denotes the multiplier ideal sheaf of the pulled back metric
$(\pi^\ast h)^m$
on
$\tilde{C}$.
\end{Def}

\noindent
This definition leads directly to the birational invariance of the 
intersection numbers, i.e.\/for a birational morphism 
$f: \widetilde{X} \rightarrow X$
\[ (f^\ast L,h^\ast).\bar{C} = (L,h).C \]
where
$\bar{C}$
is the strict transform (s. subsection~\ref{birinv-subsec}).

\noindent
Next, one has
\begin{prop} \label{Inteq-prop}
If 
$C$
is smooth,
\[ (L,h).C = L.C - \sum_{x \in C} \nu(\Theta_{h|C}, x), \]
where
$\nu(\Theta_{h|C}, x)$
is the Lelong number of the positive current 
$\Theta_h$
restricted to
$C$
in
$x \in C$.
\end{prop}

\noindent
This equality gives a more geometric interpretation of the intersection 
numbers (especially in the case of analytic singularities, s. 
Proposition~\ref{analsing-prop}) and is an important
step towards a last equality. This is the most subtle one, and to formulate 
it properly, one has to remember that plurisubharmonic functions are equal to
$-\infty$
on pluripolar sets, and that for a positive current 
$\Theta$,
the level sets of the Lelong numbers 
$E_c(\Theta) = \left\{ x \in X | \nu(\Theta,x) \geq c \right\}$
are analytic subsets of
$X$ (\cite{Siu74},\cite[(2.10)]{Dem00}):
\begin{Def} \label{Lhgen-def}
Let
$X,L,h$
be as in the previous definition. A smooth curve
$C \subset X$
will be called 
$(L,h)-\!$
general iff
$h_{|C}$
is a well defined singular metric on
$C$
and
\begin{itemize}
\item[(i)]
$C$ 
intersects no codim-2-component in any of the 
$E_c(h)$,
\item[(ii)]
$C$ 
intersects every prime divisor
$D \subset E_c(h)$
in the regular locus 
$D_{\mathrm{reg}}$
of this divisor, 
$C$
does not intersect the intersection of two such prime divisors, and every 
intersection point 
$x$
has the 
minimal Lelong number
$\nu(h,x) = \nu(h,D) := \min_{z \in D} \nu(h,z)$,
\item[(iii)]
for all
$x \in C$,
the Lelong numbers
\[ \nu(h_{|C},x) = \nu(h,x). \] 
\end{itemize} 
\end{Def}

\noindent
Using methods of \cite{BMElM00}
it is possible to show that in families of curves covering
$X$
(e.g. the Chow variety) 
every curve outside a pluripolar set is
$(L,h)-\!$
general, s. Theorem~\ref{FamLh-thm}. One can even proof the stability of this 
notion under certain blow 
ups, s. Lemma~\ref{birinvLh-lem}. The main reason for introducing this notion 
lies in the equality
\[ \mathcal{I}(h^m) \cdot \mathcal{O}_C =  \mathcal{I}(h^m)_{|C} = 
   \mathcal{I}(h^m_{|C}), \]
which is true for
$(L,h)-\!$
general curves.
From this one easily gets the announced last equality
\begin{thm} \label{Inteq-thm}
For
$(L,h)-\!$
general smooth curves
$C \subset X$,
\[ (L,h).C = \limsup_{m \rightarrow \infty} \frac{1}{m} 
              h^0(C, \mathcal{O}_C(mL) \otimes 
                             \mathcal{I}(h^m) \cdot \mathcal{O}_C), \]
where 
$\mathcal{I}(h^m) \cdot \mathcal{O}_C$
is the image of
$\mathcal{I}(h^m) \otimes \mathcal{O}_C$
in
$\mathcal{O}_C$.
\end{thm}

\noindent
This equality is 
needed in order to be able to interchange restriction (to curves
$C$)
with taking global sections (of the sheaf
$\mathcal{O}_X(mL) \otimes \mathcal{I}(h^m)$) as in the proof of the Key Lemma
in subsection~\ref{keylem-subsec}.
There are explicit counterexamples for arbitrary curves, s. 
subsection~\ref{countex-ssec}. On the other hand the
equality is true in general in case of analytic singularities, s. 
Proposition~\ref{Inteq-anal-prop}.

\noindent
Furthermore this equality motivates the introduction of intersection numbers 
as defined above, because it associates these numbers with the sheaf
$\mathcal{O}_X(mL) \otimes \mathcal{I}(h^m)$
which is useful and interesting in several ways: For certain L's and h's it 
occurs in Nadel's Vanishing Theorem. Next, an appropriately chosen metric 
$h$
should encode the ``negative part'' of 
$L$
in the multiplier ideal sheaves
$\mathcal{I}(h^m)$. 
And for an analytic Zariski decomposition 
$h$
the global sections of
$\mathcal{O}_X(mL) \otimes \mathcal{I}(h^m)$
are the same as those of
$\mathcal{O}_X(mL)$. 
 
\noindent
Second, this note proves the existence of the reduction map with 
respect to the pair
$(L,h)$. The aim is to get a reduction map with numerically trivial fibers:
\begin{Def}
Let
$X$
be a smooth projective complex manifold, let
$L$
be a pseudoeffective holomorphic line bundle on
$X$
with positive singular hermitian metric
$h$. 
Then a subvariety
$Y \subset X$
is called numerically trivial (with respect to 
$(L,h)$)
if each curve
$C \subset Y$
such that 
$h_{|C} \not\equiv \infty$
has intersection number
$(L,h).C = 0$.
\end{Def}

\noindent
If
$X$
itself is numerically trivial one can prove the following consequence for the 
curvature current:
\begin{thm} \label{NumTriv-thm}
Let
$X$
be a smooth projective complex manifold, let
$L$
be a pseudo-effective line bundle on
$X$
with positive singular hermitian metric 
$h$
such that
$X$
is
$(L,h)-\!$
numerically trivial. Then the curvature current
$\Theta_h$
may be decomposed as
\[ \Theta_h = \sum_i a_i [D_i] \]
where the
$D_i$
form a countable set of prime divisors on
$X$ 
and the 
$a_i > 0$. 
\end{thm}

\noindent
One has to adjust Tsuji's statement of the Reduction Map Theorem:
\begin{thm}
Let
$X$
be a smooth projective complex manifold, let
$L$
be a pseudoeffective holomorphic line bundle on
$X$
with positive singular hermitian metric
$h$.
Then there exists a dominant rational map
$f: X \dasharrow Y$
with connected fibres such that
\begin{itemize}
\item[(i)]
$(L,h)$ 
is numerically trivial on fibres over points of
$Y$
outside a pluripolar set. 
\item[(ii)]
for all
$x \in X$
outside a pluripolar set, every curve
$C$
through
$x$ 
with 
$\dim f(C) > 0$
has intersection number
$(L,h).C > 0$ 
\end{itemize}
Here, fibres of 
$f$
are fibres of the graph
$\Gamma_f \subset X \times Y \rightarrow Y$
seen as subschemes of 
$X$.

\noindent
Finally,
$f$ 
is uniquely determined up to birational equivalence of 
$Y$.
\end{thm}

\noindent
There are two main differences to \cite{Tsu00}. First, the reduction map need 
not be almost holomorphic. A counter example was already given in 
\cite{BCEKPRSW00}. Second, Tsuji ignores the fact that the singularities of 
arbitrary positive singular hermitian metrics are even more complicated than 
very general algebraic sets: they lie on pluripolar sets. This means for 
example, that the restriction of the singular metric may be well defined only
on fibres over points outside a pluripolar set. But this is not
so bad: for example, the Zariski closure of the union of these 
fibres is always the whole variety. 

\noindent
After these adjustments it is possible to apply Tsuji's ideas in proving the 
Reduction Map Theorem:
\begin{itemize}
\item[(a)]
For each ample divisor
$H$
and each pair 
$(L,h)$
of a line bundle with a positive singular hermitian metric one can define a
\textbf{volume}
\[ \mu_h(X, H+mL) := (\dim X)! \limsup_{l \rightarrow \infty} l^{-\dim X}
                 h^0(X, \mathcal{O}_X(l(H+mL)) \otimes \mathcal{I}(h^{ml})) \]
and we have the following

\noindent
\begin{lem} \label{vol-lem}
\[ (L,h) \mathrm{\ numerically\ not\ trivial\ \ } 
   \Rightarrow \ \limsup_{m \rightarrow \infty} \mu_h(X, H+mL) = \infty. \] 

\end{lem}
\item[(b)]
The lemma implies that for all 
$N$
there exists an
$m_0$
such that for arbitrarily large
$l \gg 0$
there exist sections
\[ 0 \not\equiv \sigma_l \in H^0(X, \mathcal{O}_X(l(H+m_0L)) \otimes 
                                    \mathcal{I}(h^{m_0l}) \otimes 
                                    \hbox{\gothic{m}}_x^{Nl}) \]
for a sufficiently general point
$x \in X$.

\item[(c)]
This is used for
\begin{keylem} \label{Key-lem}
Let
$f: M \rightarrow B$
be a projective surjective morphism from a smooth variety
$M$
to a smooth curve 
$B$. 
Let
$(L,h)$
be a pseudo-effective line bundle 
$L$
with positive singular hermitian metric
$h$.
Suppose that
$(L,h)$
is numerically trivial on all fibres
$F$
of
$f$ 
over a set 
$B^\prime \subset B$
not of Lebesgue measure 
$0$. 
If furthermore there is an
$(L,h)-\!$
general curve 
$W$
with
$f(W) = B$,
$(L,h)$
is numerically trivial on
$W$,
then
$(L,h)$
will be numerically trivial on
$M$.
\end{keylem}

\noindent
The proof is done by contradiction: Any
$\sigma_l$
as above must be
$0$.

\item[(d)]
Finally the theorem is derived from the Key Lemma with methods similar to 
those in
\cite{BCEKPRSW00}.
\end{itemize}
The intersection number equality in Theorem~\ref{Inteq-thm} is needed
essentially in proving the Key Lemma~\ref{Key-lem}, while the definition of 
the intersection number is used several times for switching to birationally
equivalent varieties.

\subsubsection*{Acknowledgement }

This note owes its subject and many ideas to the preprint \cite{Tsu00} of 
H. Tsuji. He also kindly answered some questions of the author during the
conference of the DFG Research center ``Global Methods in Complex Geometry''
in Marburg, June 2001.
This DFG research center also supported a small workshop in January 2001, where
the author had its first contact with \cite{Tsu00}. Finally, he would like to 
thank Jean-Pierre Demailly who taught him everything 
about the analytic side of multiplier
ideal sheaves during the ICTP summer school ``Vanishing theorems and effective
results in algebraic geometry'' in April 2000 and whose wonderful lecture 
notes really served, as intended, as an ``analytic tool box''. Furthermore, he 
pointed out to the author a serious mistake in a previous version of this paper
and showed him the paper \cite{BMElM00}.
Similarly, Robert Lazarsfeld on the same summer school introduced the author
to algebraic multiplier ideals.

\subsubsection*{Notation}

\noindent
If not otherwise stated, in the following
$X$
is a smooth projective complex manifold, 
$L$
a holomorphic line bundle on
$X$
with positive singular hermitian metric
$h$. 
On any trivialization
$\theta: L_{|\Omega} \stackrel{\simeq}{\rightarrow} \Omega \times \mathbb{C}$,
this metric is given by
\[ \| \xi \| = | \theta(\xi) | e^{-\phi(x)},\ \ x \in \Omega, \xi \in L_x. \]
The function
$\phi \in L^1_{\mathrm{loc}}(\Omega)$ 
is called the \textbf{weight} 
of the metric with respect to the trivialization
$\theta$,
and
$h$
is positive iff 
$\phi$
is plurisubharmonic. In that case, the \textbf{Lelong number} of 
$h$
(or 
$\phi$)
in
$z_0 \in X$
is defined as
\[ \nu(h,z_0) = \nu(\phi,z_0) = \liminf_{z \rightarrow z_0} 
                \frac{\phi(z)}{\log |z-z_0|}. \] 
If
$C \subset X$
is a smooth curve and
$h_{|C} \not\equiv \infty$,
the metric
$h_{|C}$
will be a positive singular hermitian metric on
$C$,
and there are Lelong numbers 
$\nu(h_{|C},x)$
for points
$x \in C$.

\section{Intersection numbers}

\noindent
The aim of this section is to prove Proposition~\ref{Inteq-prop}, 
Theorem~\ref{NumTriv-thm} and 
Theorem~\ref{Inteq-thm}. To this purpose one has to study the behavior of 
slices of positive currents, especially what happens to their
Lelong numbers.

\noindent
Furthermore, the 
birational invariance of intersection numbers is shown, and  
intersection numbers are computed on log resolutions
in case of analytic singularities.

\subsection{Slices of positive currents} \label{slices-ssec}

\noindent
The aim is to prove the following
\begin{thm} \label{FamLh-thm}
Let
$\pi: X \rightarrow B$
be a smooth family
$X$
of smooth projective curves over a smooth quasiprojective base
$B$.
Let
$L$
be a pseudo-effective line bundle on 
$X$
and
$h$ 
a positive singular hermitian metric on
$L$.
Then there is a pluripolar set
$N_B \subset B$
such that for
$a \in B - N_B$,
every fibre
$\pi^{-1}(a)$,
is an
$(L,h)-\!$
general curve.
\end{thm}

\noindent
This Theorem is essentially a consequence of Ben Messaoud's
\begin{thm} \label{Tranche-thm}
Let
$M$,
$G$
be two complex varieties of dimension
$n$
and
$k$,
let 
$\phi$
be a plurisubharmonic function on
$M$
and let
$f: M \rightarrow G$
be a submersion admitting a holomorphic section
$s$.
Then there exists a pluripolar set
$E \subset G$
such that for all
$a \in G \setminus E$,
the restricted pluripolar function
$\phi_{|f^{-1}(a)} \not\equiv -\infty$
and
\[ \nu(\phi, s(a)) = \nu(\phi_{|f^{-1}(a)}, s(a)). \]
\end{thm}
\begin{proof}
S. \cite[Cor. 5.4]{BMElM00}.
\end{proof}

\noindent
Now take an open subset
$U \subset X$
such that
$\pi: U \cong \Delta^k \times \Delta \rightarrow \Delta^k$
with
$\Delta \subset \mathbb{C}$
the unit disk. Apply Theorem~\ref{Tranche-thm} to the family
$U \times \Delta \stackrel{\pi \times \mathrm{id}_\Delta}{\longrightarrow} 
 \Delta^k \times \Delta$,
the pulled back plurisubharmonic function and the section
\[ s : \Delta^k \times \Delta \rightarrow U \times \Delta,\ 
       (b,t) \mapsto (b,t,t). \]
Since the projection of
$s(\Delta^k \times \Delta)$
on
$U$
is an isomorphism there is pluripolar set
$E_U \subset U$
such that for all
$x \in U \setminus E_U$
\[ \nu(\phi,x) = \nu(\phi_{|\pi^{-1}(\pi(x))},x). \]
Since the countable union of pluripolar sets is again pluripolar the same is 
true for a pluripolar set
$E \subset X$.
The other two requirements of Definition~\ref{Lhgen-def} for an
$(L,h)-\!$
general curve also show that these curves must be fibres outside the countable
union of analytic subsets, which is a pluripolar set. This shows the Theorem.

\subsection{Proof of Proposition~\ref{Inteq-prop}} 
\label{first-eq-ssec}

The first step is to compare the sum of the restricted Lelong numbers on
arbitray curves
$C \in X$
with
$h_{|C} \not\equiv  \infty$
to the ordinary intersection number of
$C$
with
$L$:
\begin{lem} \label{Lel=intnum-lem}
\[ \sum_{x \in C} \nu(h_{|C},x) \leq L.C . \]
\end{lem}
\begin{proof}
Since 
$h_{|C}$
is positive, the curvature current
$i\Theta_{h_{|C}} \geq 0$,
too.  By a theorem of Siu, the
Lelong level sets 
$E_c(\phi) = \left\{ x \in X: \nu(\phi,x) \geq c \right\}$
are analytic \cite[(2.10)]{Dem00}. But then there are only countably many 
points 
$(x_i)_{i \in \mathbb{N}}$
on 
$C$
with
$\nu(h_{|C},x_i) \neq 0$.
By Siu's decomposition formula \cite[(2.18)]{Dem00}
the current
$i\Theta_{h_{|C}} - \sum_{i=1}^N \nu(h_{|C},x_i)[x_i]$
is still positive for arbitrary
$N$
(where
$[x_i]$ is the integration current of the divisor
$x_i$).
Consequently the first Chern class of the 
($\mathbb{R}-\!$)
divisor
$L_{|C} - \sum_{i=1}^N \nu(h_{|C},x_i)x_i$
is
$\geq 0$,
hence 
$L.C - \sum_{i=1}^N \nu(h_{|C},x_i) \geq 0$,
and the claim follows. 
\end{proof}

\begin{lem}
Let
$C$
be a smooth curve and
$h$
a positive singular hermitian metric on
$C$.
Then:
\[ \limsup_{m \rightarrow \infty} \frac{1}{m} \deg_C \mathcal{I}(h^m) =
   - \sum_{x \in C} \nu(h,x). \]
\end{lem}
\begin{proof}
$\mathcal{I}(h^m)$
is a torsion free subsheaf of
$\mathcal{O}_C$,
hence it corresponds to a divisor on
$C$,
say
$\mathcal{I}(h^m) = \mathcal{O}(-D_m)$,
where
$D_m$
is an effective divisor on
$C$.
We show that
\begin{eqnarray} 
   & \mathrm{mult}_x D_m \leq \nu(h^m,x) <  
   \mathrm{mult}_x D_m + 1 & .
   \label{m=m+1-eq}
\end{eqnarray}
This is true for arbitrary positive metrics 
$h$:
Choose a sufficiently small neighborhood
$U$
of 
$x$
such that
$\nu(h,y) < 1$
for all 
$y \in U\setminus\{x\}$.
Let
$\phi_h$
and
$\Theta_h$
be the plurisubharmonic function and
$(1,1)-\!$
current corresponding to
$h$
in
$U$. 
As explained in the proof of the previous lemma the current
$\Theta = \Theta_h - \nu(h,x)[x]$
is still positive, with
$\nu(\Theta,x) = 0$,
$\nu(\Theta,y) < 1$
for all
$y \in U\setminus\{x\}$.
Let
$\psi$
be a plurisubharmonic function with
$dd^c\psi = \Theta$.
Then
$\phi_h = \psi + \nu(h,x)\log|z-x|$,
hence
$-\mathrm{mult}_x \mathcal{I}(h) \geq \lfloor \nu(h,x) \rfloor$.

\noindent
On the other hand
$e^{-2(\psi + (\nu(h,x)-\lfloor \nu(h,x) \rfloor)\log|z-x|)}$ 
is locally integrable around 
$x$
since the Lelong number in
$x$
is
$<1$,
by Skoda's lemma~\cite[(5.6)]{Dem00}. This proves
$\mathcal{I}(h)_x = ((z-x)^{\lfloor \nu(h,x) \rfloor})$, 
hence (\ref{m=m+1-eq}).

\noindent
Now one concludes:
\begin{eqnarray*} 
   \liminf_{m \rightarrow \infty} \frac{1}{m} \deg_C D_m 
   \stackrel{(\ref{m=m+1-eq})}{=} 
   \liminf_{m \rightarrow \infty} \frac{1}{m} \sum_{x \in C}
   \left\lfloor \nu(h^m,x) \right\rfloor & \leq &
   \liminf_{m \rightarrow \infty} \frac{1}{m} \sum_{x \in C} \nu(h^m,x) \\
     & = &  \sum_{x \in C} \nu(h,x) < \infty.
\end{eqnarray*}
The inequality
$a - \frac{1}{m} \left\lfloor ma \right\rfloor < \frac{1}{m}$
for arbitrary
$a \in \mathbb{R}, m \in \mathbb{N}$
shows the lemma.
\end{proof}

\noindent
Proposition~\ref{Inteq-prop} follows from
\begin{eqnarray}
\ \ \ \ \ \limsup_{m \rightarrow \infty} \frac{1}{m} h^0(C, \mathcal{O}_C(mL) 
                                             \otimes \mathcal{O}_C(-D_m)) 
& = & \limsup_{m \rightarrow \infty} \frac{1}{m} \deg_C (mL - D_m) 
\label{h0=deg-new-eq} 
\end{eqnarray}
\begin{proof}
By Lemma~\ref{Lel=intnum-lem},
$\deg_C D_m \leq \sum_{x \in C} \nu(h_{|C}^m, x) \leq mL.C$. 
Consequently,
\[ \deg_C (mL - D_m) \geq 0. \] 
Let
$g(C)$
be the genus of the curve
$C$. If
$\deg_C (mL - D_m) \leq 2g(C) - 2$
and 
$mL - D_m$
is not effective,
$H^0(C, \mathcal{O}(mL - D_m)) = 0$.
If
$\deg_C (mL - D_m) \leq 2g(C) - 2$
and 
$mL - D_m$
is effective, 
\[ H^0(C, \mathcal{O}(mL - D_m)) \leq \deg_C (mL - D_m) + 1 \leq 2g(C) - 1. \]
If 
$\deg_C (mL - D_m) > 2g(C) - 2$,
then
$H^1(C, \mathcal{O}(mL - D_m)) = 0$,
and (\ref{h0=deg-new-eq}) will follow by Riemann-Roch. 
\end{proof}

\subsection{Proof of Theorem~\ref{Inteq-thm}}

One main ingredient of the proof, which is useful in many circumstances, is 
\begin{Ext-thm}[Ohsawa-Takegoshi]
Let
$\Omega \subset \mathbb{C}^n$ 
be a bounded open pseudoconvex set,
$L = \{ z_i = \ldots z_n = 0\}$,
$1 \leq i \leq n$,
a linear subspace, and
$\psi \in \mathrm{Psh}(\Omega)$
with
$\psi_{|L} \neq -\infty$.

\noindent
Then there is a constant
$C > 0$,
only depending on
$n$,
such that for all 
holomorphic functions 
$f$
on
$L$
with
$\int_{L \cap \Omega} |f|^2 e^{-2\psi}d\lambda_L < \infty$,
there is an
$F \in \mathcal{O}(\Omega)$
such that
$F_{|L} = f$
and
\[ \int_\Omega |F|^2 e^{-2\psi}d\lambda_\Omega \leq
   C \cdot \int_{L \cap \Omega} |f|^2 e^{-2\psi}d\lambda_L \]
\end{Ext-thm}
\begin{proof}
S. \cite[(12.9)]{Dem00}.
\end{proof}

\noindent
Now let
$C$
be a smooth 
$(L,h)-\!$
general fibre curve of the family
$X$.
Let 
$D_m, D_m^\prime$
be the effective divisors corresponding to the ideal sheaves
$\mathcal{I}(h_{|C}^m)$
and
$\mathcal{I}(h^m)_{|C}$,
as explained in subsection~\ref{first-eq-ssec}.
The Extension Theorem implies a natural inclusion
\[ \mathcal{I}(h_{|C}^m) \subset \mathcal{I}(h^m)_{|C}, \]
hence
$\deg_C D_m^\prime \leq \deg_C D_m$,
and one can prove
\begin{eqnarray}
\ \ \ \ \ \limsup_{m \rightarrow \infty} \frac{1}{m} h^0(C, \mathcal{O}_C(mL) 
                                     \otimes \mathcal{O}_C(-D_m^\prime)) 
& = & \limsup_{m \rightarrow \infty} \frac{1}{m} \deg_C (mL - D_m^\prime)
\label{h0=degprime-new-eq}
\end{eqnarray}
similarly as (\ref{h0=deg-new-eq}).

\noindent
The 
$(L,h)-\!$
generality implies
\[ \mathcal{I}(h_{|C}^m) = \mathcal{I}(h^m)_{|C}. \]
\begin{proof}
By Skoda's Lemma~\cite[(5.6)]{Dem00},
$\mathcal{I}(h_{|C}^m)_x = \mathcal{I}(h^m)_{|C,x} = \mathcal{O}_{C,x}$
for all points
$x \in C$
with
$\nu(h,x) = \nu(h_{|C},x) = 0$.

\noindent
Let
$x \in C$
be a point with
$\nu = \nu(h,x) = \nu(h_{|C},x) > 0$.
Then there is a divisor
$D$
through
$x$
locally defined by
$g \in \mathcal{O}_{X,x}$,
$x \in D_{\mathrm{reg}}$
and
$\nu = \nu_{h,x} = \nu(h,D)$. 
As in subsection~\ref{slices-ssec} it follows that
\[ \mathcal{I}(h^m)_x = (g^{\left\lfloor m\nu \right\rfloor}) \subset 
   \mathcal{O}_{X,x}. \]
Similarly one concludes
\[ \mathcal{I}(h^m_{|C})_x = (g^{\left\lfloor m\nu \right\rfloor}) \subset 
   \mathcal{O}_{C,x}. \]
\end{proof}

\subsection{Characterization of numerically trivial varieties}

\noindent
The proof of Theorem~\ref{NumTriv-thm} starts with the Siu decomposition of 
the curvature current \cite[(2.18)]{Dem00}
\[ \Theta_h = \sum_i a_i [D_i] + R \]
where the 
$D_i$
are the (countably many) prime divisors in the Lelong number level sets
$E_c(h)$ 
and the
$a_i = \min_{x \in D_i} \nu(\Theta_h, x)$.

\noindent
Next, take a very ample divisor
$H$. 
By Theorem~\ref{FamLh-thm} there is a smooth complete intersection curve
$C = H_1 \cap \ldots \cap H_{n-1}$,
$H_i \in |H|$
which is
$(L,h)-\!$
general. Now by Proposition~\ref{Inteq-prop}
\[ 0 = (L,h).C = L.C - \sum_{x \in C} \nu(\sum_i a_i [D_i]_{|C},x) -
                       \sum_{x \in C} \nu(R_{|C},x).               \]
Since
$C$
is
$(L,h)-\!$
general the only points 
$x \in C$
where 
$\nu(\Theta_h,x) > 0$
are the intersection points with the regular part of the 
$D_i$'s 
where furthermore 
$\nu(\Theta_h,x) = \nu(\Theta_h,D_i) = a_i$.
Consequently
\[ 0 = (L,h).C = L.C - \sum_i a_i D_i.C. \]
But this implies 
\[ 0 = R.C = \int_X R \wedge (\omega_H)^{n-1} \]
where 
$\omega_H$
is the strictly positive
$\mathcal{C}^\infty-\!$
metric belonging to the very ample divisor
$H$.
Since
$R$
is a positive current it follows
\[ R = 0. \]

\subsection{A counterexample for non-$(L,h)-\!$ general curves}
\label{countex-ssec}

\noindent
First, one constructs a convex function
$\chi: \mathbb{R} \rightarrow \mathbb{R}$
with slow growth at
$-\infty$
(i.e. the derivation tends to 
$0$)
such that 
$\chi(-\infty) = -\infty$.
For example, take
\[ \chi(x) = \left\{ \begin{array}{ll} 
                     x & \mathrm{for\ } x \geq -1 \\
                     -\sum_{k=1}^n \frac{1}{k} + (x+n)\frac{1}{n+1} &
                     \mathrm{for\ } -n-1 \leq x \leq -n 
                     \end{array}       \right. \]
Then one considers the plurisubharmonic function 
$\psi = \max(\log|z_1|, \chi(\log|z_2|))$
on
$\mathbb{C}^2$.
The Lelong numbers
$\nu(\psi,x)$
are 
$0$
everywhere because of the slow growth of
$\chi$
at
$-\infty$,
but the restriction of
$\psi$
onto
$C = \left\{ z_2 = 0 \right \}$
has Lelong number
$\nu(\psi_{|C},0) = 1$.

\noindent
The induced metric
$h$
may be extended to a metric of the relatively ample line bundle 
$\mathcal{O}(1)$
on the 
$\mathbb{P}^1-\!$
bundle
$\mathbb{C} \times \mathbb{P}^1$
which yields the counterexample.

\subsection{Birational invariance}
\label{birinv-subsec}

Since the intersection numbers
$(L,h).C$
are computed by pulling back to the normalization
$\widehat{C}$
it is obvious that the intersection number
$(\pi^\ast L, \pi^\ast h).\bar{C}$
of the strict transform
$\bar{C}$ 
of a birational map
$\pi$
does not change. The aim of this subsection is to generalize this observation 
and to apply it in the case of analytic singularities, thus obtaining a more 
algebraic definition of the intersection numbers.
\begin{lem}
Let
$\mu: C^\prime \rightarrow C$
be a finite morphism between smooth curves. Let
$(L,h)$
be a pseudo-effective line bundle on
$C$
with
$i\Theta_h \geq 0$.
Then
\[ (\mu^\ast L, \mu^\ast h).C^\prime = \deg \mu \cdot (L,h).C. \]
\end{lem}
\begin{proof}
It is enough to consider the following situation: Let
$\mu: \Delta \rightarrow \Delta,\ z \mapsto z^n$
be a finite morphism on the unit disc 
$\Delta$,
let
$L$
be trivial on
$\Delta$
and let the function
$\psi \in \mathrm{Psh}(\Delta)$,
$i\Theta_h = dd^c \psi$.
Then
\[ \nu(\psi,0) = \liminf_{|z| \rightarrow 0} \frac{\psi(z)}{\log|z|} =
                 \liminf_{|z| \rightarrow 0} \frac{\psi(z^n)}{\log|z^n|} = 
                 \frac{1}{n} \liminf_{|z|  \rightarrow 0} 
                 \frac{\psi(z^n)}{\log|z|} =
                 \frac{1}{n} \nu(\psi \circ \mu,0). \]
Now the lemma follows by Proposition~\ref{Inteq-prop}.
\end{proof}

\begin{prop} \label{numtrivinv-prop}
Let
$f: Y \rightarrow X$
be a surjective morphism between smooth and projective varieties 
$X$
and
$Y$.
Let
$(L,h)$
be a pseudo-effective line bundle on
$X$
with
$i\Theta_h \geq 0$.
Then
\[ (L,h)\ \mathrm{numerically\ trivial\ on\ } X\ \Longleftrightarrow
   (f^\ast L,f^\ast h)\ \mathrm{numerically\ trivial\ on\ } Y. \]
\end{prop}
\begin{proof}
Assume first that 
$(L,h)$
is numerically trivial on 
$X$.
Let
$C \subset Y$
be an irreducible curve on
$Y$
with
$f^\ast h_{|C} \not\equiv \infty$.
When
$f(C)$
is a point, this point will lie in the smooth part of
$h$,
and there won't be any singularity of
$h$
on
$C$. Consequently,
\[ (f^\ast L,f^\ast h).C = f^\ast L.C = L.f_\ast C = 0. \]

\noindent
When
$f(C)$
is another irreducible curve
$C^\prime$
then one can lift the morphism
$f_{|C}$
to the smooth normalizations
$\widehat{C},\widehat{C}^\prime$,
and the above equality follows by the lemma.

\noindent
Similarly, assume that
$(f^\ast L,f^\ast h)$
is numerically trivial on
$Y$.
Let
$C$ 
be an irreducible curve on
$X$
with
$h_{|C} \not\equiv \infty$.
Then there exists an irreducible curve
$C^\prime \subset Y$
such that 
$f(C^\prime) = C$,
and the argument is as above. 
\end{proof}

\noindent
The birational invariance can be used to prove the birational invariance of
$(L,h)-\!$
generality:
\begin{lem} \label{birinvLh-lem}
Let
$C$
be a smooth
$(L,h)-\!$
general curve on 
$X$,
and let
$Z \subset X$
be a smooth subvariety with
$C \not\subset Z$,
let
$\pi: \widehat{X} \rightarrow X$
be the blowup of 
$X$
with centre 
$Z$.
Then the strict transform
$\widehat{C}$
of
$C$
is still
$(\pi^\ast L, \pi^\ast h)-\!$
general.
\end{lem} 
\begin{proof}
The assertion is clear as long as 
$Z \cap C = \emptyset$.
Otherwise, let
$x \in Z \cap C$
be a point such that
$\nu(h_{|C},x) = 0$.
Then for
$y$
the unique preimage of
$x$
in
$\widehat{C}$,
\[ 0 = \nu(\pi^\ast h, y) \leq \nu(\pi^\ast h_{|\widehat{C}},y) = 
       \nu(h_{|C},x) = 0. \]

\noindent
If
$x \in Z \cap C$
is  a point such that
$\nu(h_{|C},x) > 0$
then 
$C$
will intersect transversally a prime divisor
$D$
of some
$E_c(h)$.
Consider two cases:
\begin{itemize}
\item[(a)]
$Z$
is a point. Then the intersection of the strict transforms 
$\widehat{C} \cap \widehat{D} = \emptyset$,
and
$C$
intersects the smooth exceptional divisor
$E$
transversally in a unique point 
$y \in E$
with
$\pi(y) = x$.
Furthermore,
\[ \nu(\pi^\ast h, E) \geq \nu(h,D) = \nu(h,x) = \nu(h_{|C},x) =
   \nu(\pi^\ast h_{|\widehat{C}},y), \]
hence
$\nu(\pi^\ast h, E) = \nu(\pi^\ast h_{|\widehat{C}},y)$.
\item[(b)]
$\dim Z \geq 1$. 
Then 
$\widehat{C} \cap \widehat{D}$
consists of one point
$y$,
and by the same argument as in (a), replacing 
$E$
by
$\widehat{D}$,
it follows
\[ \nu(\pi^\ast h, \widehat{D}) = \nu(\pi^\ast h_{|\widehat{C}},y). \]
$\widehat{D}$
cannot be singular in
$y$
since then
$\nu(\pi^\ast h_{|\widehat{C}},y) \geq \nu(\pi^\ast h,y) > 
\nu(\pi^\ast h, \widehat{D})$.
\end{itemize}
\end{proof}

\subsection{The case of analytic singularities}

\noindent
The 
$(L,h)-\!$
intersection numbers are much easier to handle if the plurisubharmonic weight 
of the metric
$h$ 
has only analytic singularities:
\begin{Def}
$\phi \in \mathrm{Psh}(\Omega)$,
$\Omega \subset \mathbb{C}^n$
open, 
is said to have \textbf{analytic singularities}, if locally,
$\phi$
can be written as
\[ \phi = \frac{\alpha}{2} \log (\sum |f_i|^2) + v,\ \alpha \in \mathbb{R}^+,\]
where
$v$
is locally bounded, and the
$f_i$
are (germs of) holomorphic functions. 
\end{Def}

\noindent
For example, in this case Theorem~\ref{Inteq-thm} is true for arbitrary smooth
curves. Furthermore it is easier to  compute 
$(L,h)-\!$
intersection numbers on log resolutions.

\noindent
But first some more properties of metrics with analytic
singularities: By definition,  
the corresponding plurisubharmonic weight may locally
be written as
$\phi_h = \frac{\alpha}{2} \log (\sum |f_i|^2) + O(1)$.
Define
$\mathcal{J}(h/\alpha)$
as the ideal sheaf of germs of holomorphic functions 
$f$
such that
\[  |f| \leq C \cdot (\sum |f_i|). \] 
One can easily prove that 
$\mathcal{J}(h/\alpha)_x$
is the integral closure of the ideal generated by the germs
$f_i$
(cf. \cite[(1.11)]{Dem00}).
Consequently, 
$\mathcal{J}(h/\alpha)$
is coherent. Furthermore,
\[ \mathcal{J}(h/\alpha)_x = (g_1, \ldots, g_M) \Longrightarrow 
   \phi = \frac{\alpha}{2} \log (\sum |g_i|^2) + O(1). \]
There exists log resolutions 
$\mu: X^\prime \rightarrow X$
of
$\mathcal{J}(h/\alpha)$
with
$X^\prime$
non-singular, i.e.
\begin{itemize}
\item[(a)]
$\mu$
is proper birational,
\item[(b)]
$\mu^{-1}\mathcal{J}(h/\alpha) = 
 \mathcal{J}(h/\alpha) \cdot \mathcal{O}_{X^\prime} = 
 \mathcal{O}_{X^\prime}(-F)$
where 
$F$
is an effective divisor on
$X^\prime$
such that
$F+\mathrm{Exc}(\mu)$
has simple normal crossing support.
\end{itemize} 

\noindent
An existence proof is contained in the Hironaka package, cf.~\cite{BM97}. 

\noindent
The main tool when dealing with metrics with analytic singularities is the 
following theorem which may be seen as an algebraic definition of multiplier 
ideals:
\begin{thm}
$\mathcal{I}(h) = \mu_\ast(K_{\bar{\Omega}/\Omega} - [\alpha F])$.
\end{thm}
\begin{proof}
See~\cite[(5.9)]{Dem00}.
\end{proof}

\noindent
The aim is now to prove
\begin{prop} \label{Inteq-anal-prop}
Let 
$X$
be a quasi-projective manifold,
$L$
a pseudo-effective line bundle with positive singular hermitian metric
$h$.
Then for every smooth curve
$C \subset X$,
\[ (L,h).C = \limsup_{m \rightarrow \infty} \frac{1}{m} 
              h^0(C, \mathcal{O}_C(mL) \otimes 
                             \mathcal{I}(h^m) \cdot \mathcal{O}_C), \]
where 
$\mathcal{I}(h^m) \cdot \mathcal{O}_C$
is the image of
$\mathcal{I}(h^m) \otimes \mathcal{O}_C$
in
$\mathcal{O}_C$.
\end{prop}
\begin{proof}
Let
$D_m, D_m^\prime$
be effective divisors corresponding to the torsion free ideal sheaves
$\mathcal{I}(h^m_{|C}), 
 \mathcal{I}(h^m) \cdot \mathcal{O}_C = \mathcal{I}(h^m)_{|C}$. 
By (\ref{h0=deg-new-eq}),(\ref{h0=degprime-new-eq}) it is enough to show that
\begin{eqnarray} 
\liminf_{m \rightarrow \infty} \frac{1}{m} \deg_C D_m & = &
\liminf_{m \rightarrow \infty} \frac{1}{m} \deg_C D_m^\prime
\label{deg-eq}
\end{eqnarray}

\noindent
Let
$(x_i)_{i \in \mathbb{N}}$
be the countably many points on
$C$
such that
$\mathrm{mult}_x D_m \neq 0$
or
$\mathrm{mult}_x D_m^\prime \neq 0$
for some
$m \in \mathbb{N}$.
Since
$C$
is smooth there is an open subset
$U \subset X$
containing all the 
$x_i$
such that
$C = H_2 \cap \ldots \cap H_n$
is a complete intersection of very ample smooth hypersurfaces
$H_i \subset U$.
It is enough to prove (\ref{deg-eq}) on
$U$,
hence one can assume without loss of generality that
$C$
is a complete intersection on
$X$
as above.

\noindent
Now construct a log resolution
$\mu: X^\prime \rightarrow X$
as above such that furthermore,
\begin{itemize}
\item[(c)]
the support of 
$F$
contains the support of
$\mathrm{Exc}(\mu)$,
\item[(d)]
the strict transforms 
$H_i^\prime$
of the
$H_i$
are smooth,
$\sum H_i + F$
has simple normal crossing support and
$\mu^\ast H_n = H_n^\prime + \sum b_j E_j$
where the
$E_j$
are prime components of
$\mathrm{Exc}(\mu)$.
\end{itemize}

\noindent
One has
\begin{thm}[Local vanishing]
Let
$\hbox{\gothic{a}} \subset \mathcal{O}_X$
be an ideal sheaf on a smooth quasiprojective complex variety
$X$,
and let
$\mu: X^\prime \rightarrow X$
be a log resolution of 
$\hbox{\gothic{a}}$
with
$\hbox{\gothic{a}} \cdot \mathcal{O}_{X^\prime} = 
 \mathcal{O}_{X^\prime}(-F)$. 
Then for any rational
$c > 0$:
\[ R^j\mu_\ast \mathcal{O}_{X^\prime}(K_{X^\prime/X} - [c \cdot F]) = 0\ 
   \mathrm{for\ } j > 0. \] 
\end{thm}
\begin{proof}
See~\cite[4.3]{Laz00}.
\end{proof}

\noindent
This theorem is used to prove the following inclusions of ideal sheaves on
$H_n$:
There is a
$c \in \mathbb{N}$
independent of 
$m$
such that
\[ \mathcal{I}(h^{m+c})_{|H_n} \subset \mathcal{I}(h^m_{|H_n}) \subset 
   \mathcal{I}(h^m)_{|H_n}. \]
\noindent
Equation (\ref{deg-eq}) follows by induction and 
$\limsup_{m \rightarrow \infty} \frac{1}{m} d_{m+c} = 
 \limsup_{m \rightarrow \infty} \frac{1}{m} d_m$
for sequences 
$(d_m)_{m \in \mathbb{N}}$. 

\noindent
The proof of this fact is modelled on the proof of the Restriction Theorem 
\cite[(5.1)]{Laz00}. First of all,
$\mu_{|H_n}: H_n^\prime \rightarrow H_n$
is a log resolution of
$\mathcal{J}(h/\alpha)_{|H_n} = \mathcal{J}(h/\alpha) \cdot \mathcal{O}_{H_n}$
by property (d) of
$\mu$.
Property 
$(c)$
implies that there exists
$c \in \mathbb{N}$
such that 
\[ K_{X^\prime/X} -[m\alpha F] - c\alpha F \subset 
   K_{X^\prime/X} -[m\alpha F] - \sum b_j E_j =: B, \]
and consequently 
\[ \mathcal{I}(h^{m+c})_{|H_n} = \mu_\ast(K_{X^\prime/X} 
                                 -[(m+c)\alpha F])_{|H_n} \subset 
   \mu_\ast \mathcal{O}_{X^\prime}(B)_{|H_n}. \]
Now,
$B - H_n^\prime = K_{X^\prime/X} -[m\alpha F] - \mu^\ast H_n$.
Local vanishing applied on 
$\mathcal{J}(h/m\alpha) \cdot \mathcal{O}(-H_n)$
implies 
\[ R^1 \mu_\ast \mathcal{O}_{X^\prime}(B - H_n^\prime) = 0. \]
Then
\[ \mu_\ast \mathcal{O}_{X^\prime}(B)_{|H_n} = 
   (\mu_{|H_n})_\ast(\mathcal{O}_{H_n^\prime}(B_{|H_n^\prime})) \]
follows by taking direct images in the exact sequence
\[ 0 \rightarrow \mathcal{O}_{X^\prime}(B - H_n^\prime) 
     \stackrel{\cdot H_n^\prime}{\rightarrow}
     \mathcal{O}_{X^\prime}(B) \rightarrow 
     \mathcal{O}_{H_n^\prime}(B_{|H_n^\prime}) \rightarrow 0. \]
Since
$K_{H_n^\prime/H_n} = (K_{X^\prime/X} - \sum b_j D_j)_{|H_n^\prime}$,
it follows 
\[ (\mu_{|H_n})_\ast(\mathcal{O}_{H_n^\prime}(B_{|H_n^\prime})) =
   (\mu_{|H_n})_\ast(K_{H_n^\prime/H_n} - [m\alpha F_{|H_n^\prime}]) =
   \mathcal{I}(h^m_{|H_n}), \]
hence the first inclusion.

\noindent
The second inclusion follows by the Ohsawa-Takegoshi Extension Theorem.
\end{proof}

\noindent
The last part of this subsection shows how to compute the 
$(L,h)-\!$
intersection numbers for metrics with analytic singularities on a log 
resolution of the ideal sheaf of the singularities:
\begin{prop} \label{analsing-prop}
Let
$X$
be a smooth projective variety, let
$(L,h)$
be a pseudo-effective line bundle 
$L$
on
$X$
with a singular hermitian metric 
$h$
such that
$i\Theta_h \geq 0$
and
$h$
has analytic singularities. Let
$\mathcal{J}(h/\alpha)$
be the ideal sheaf ot these singularities, let
$\mu: \tilde{X} \rightarrow X$
be a log resolution of 
$X$
with
$\mu^\ast \mathcal{J}(h/\alpha) = \mathcal{O}(-F)$.
Let
$C \subset X$
be an irreducible curve.
Then
\[ (L,h).C = \mu^\ast L.\bar{C} - F.\bar{C}, \]
where
$\bar{C}$
is the strict transform of
$C$.
\end{prop}
\begin{proof}
By birational invariance, 
\[ (L,h).C = (\mu^\ast L,\mu^\ast h).\bar{C}. \]
But the pull back of 
$h$
is just the metric given by
$F$
by definition of analytic singularities and log resolutions
(s. \cite[(3.13)]{Dem00}). This implies the
proposition.
\end{proof}

\section{The Reduction Map Theorem}

\subsection{The volume $\mu_h$ and numerical triviality} \label{vol-ssec}

\noindent
The aim of this subsection is the prove of Lemma~\ref{vol-lem} and the 
existence of a section
$\sigma_\lambda$ 
as in step (b) of the introduction.

\noindent
The proof is by induction on 
$\dim X$.
If
$X = C$
is a smooth curve, the volume will be
\begin{eqnarray*}
\mu_h(C, H+mL) & = & \limsup_{l \rightarrow \infty} \frac{1}{l}
H^0(C, \mathcal{O}_C(l(H+mL)) \otimes \mathcal{I}(h^{ml})) = \\
 & = & \limsup_{l \rightarrow \infty} \frac{1}{l}
\deg_C \mathcal{O}_C(l(H+mL)) \otimes \mathcal{I}(h^{ml})) = \\
 & = & \deg_C H + (L^{\otimes m},h^{\otimes m}).C = \deg_C H + m\cdot(L,h).C,
\end{eqnarray*}
where the second and the third equality follow by 
equation~(\ref{h0=deg-new-eq}),
while the fourth is a consequence of the Lelong number formula for the 
intersection number.

\noindent
If
$\dim X = n$,
then for every
$n_1 \gg 0$
there will be a hyperplane pencil in
$|n_1H|$
with smooth center
$Z \subset X$
such that the general element
$F$
of the pencil is smooth, and for sufficiently general
$F$,
the restricted metric
$h_{|F} \not \equiv \infty$.

\noindent
\textbf{Step 1.}
$(L,h)$
is not numerically trivial on a sufficiently general
$F$.

\noindent
Let
$C \in X$
be an irreducible, not necessarily smooth curve such that
$(L,h).C > 0$.

\noindent
\textit{Claim.} For arbitrary 
$n_i \gg 0$,
there exists a complete intersection
\[ H_1 \cap \ldots \cap H_{n-1} = C \cup \bigcup_i C_i,\ \  H_i \in |n_i H|, \]
such that the
$C_i$
are irreducible smooth curves with
$h_{|C_i} \not\equiv \infty$.

\begin{proof}
If
$n=2$,
the curve 
$C$
is a divisor, and for
$m \gg 0$,
the linear system
$|mH-C|$
is very ample. Hence a general element
$C^\prime \in |mH-C|$
is irreducible, and
$h_{|C^\prime} \not\equiv \infty$.

\noindent
For
$n > 2$,
the curve
$C$
is contained in irreducible hypersurface 
$H^\prime$
with
$h_{|H^\prime} \not\equiv \infty$.
For some
$m \gg 0$
the linear system
$|mH-H^\prime|$
is very ample. Hence a general element
$H^{\prime\prime} \in |mH-H^\prime|$
is irreducible, and
$h_{|H^{\prime\prime}} \not\equiv \infty$.
Use induction on
$H_1 = H^\prime \cup H^{\prime\prime}$.
\end{proof}

\noindent
\textit{Claim.} For every irreducible curve
$C \subset X$, 
the following inequality is true:
\begin{eqnarray} \label{nongenLel-eq}
(L,h).C \leq L.C - \sum_j \nu(\Theta_h,D_j) C.D_j. 
\end{eqnarray}

\noindent
\begin{proof}
Let
$\pi: \widehat{C} \rightarrow C$
be the normalization of
$C$. By the decomposition theorem of Siu \cite[(2.18)]{Dem00},
\[ i\Theta_h = \sum_j \nu(\Theta_h, D_j)[D_j] + R,\ R \geq 0, \]
where 
$R$
is a positive residual 
$(1,1)-\!$
current. Let 
$\phi_j, \phi_R$
be the plurisubharmonic functions such that (locally)
$dd^c \phi_j = [D_j]$
and
$dd^c \phi_R = R$.
Then
\[ \Theta_{\pi^\ast h} = \sum_j \nu(\Theta_h, D_j) dd^c(\phi_j \circ \pi) +
                                                    dd^c(\phi_R \circ \pi) 
                    \geq \sum_j \nu(\Theta_h, D_j)[\pi^\ast D_j]. \]
Since
$\liminf_{z \rightarrow x} \sum f_j(z) \geq \sum_j \liminf_{z \rightarrow x} 
                                                                      f_j(z)$
for arbitrary functions
$f_j$,
it follows
\[ \nu(\Theta_{\pi^\ast h},x) \geq \sum_j \nu(\Theta_h, D_j)
   \nu([\pi^\ast D_j],x) = \sum_j \nu(\Theta_h, D_j) \cdot 
   \mathrm{mult}_x \pi^\ast D_j\ \ \forall x \in \widehat{C}. \]
But then
\begin{eqnarray*} 
(L,h).C = L.C - \sum_{x \in \widehat{C}} \nu(\Theta_{\pi^\ast h},x) &
 \leq & L.C - \sum_{x \in \widehat{C}} 
        (\sum_j \nu(\Theta_h, D_j) \cdot \mathrm{mult}_x \pi^\ast D_j) = \\
 & =  & L.C - \sum_j \nu(\Theta_h,D_j) C.D_j. 
\end{eqnarray*}
\end{proof}

\noindent
By Theorem~\ref{FamLh-thm}, on a sufficiently general fibre
$F$,
there is a smooth irreducible curve
\[ C_F = H_1^\prime \cap \ldots H_{n-1}^\prime,\ \ H_i^\prime \in |n_i H|, \]
which is
$(L,h)-\!$
general. Consequently,
\[ (L,h).C_F = L.C_F - \sum_j \nu(\Theta_h,D_j) C_F.D_j. \]
But this implies together with (\ref{nongenLel-eq})
\[ (L,h).C_F \geq (L,h).(C + \sum C_i) > 0, \]
because
$L.C_F = L.(C + \sum C_i)$.
Hence
$F$
is not numerically trivial.

\noindent
\textbf{Step 2.} Let
$\mu: \widetilde{X} \rightarrow X$
be the blow up of
$X$
in
$Z$,
with exceptional divisor
$E$.
Then
\[ \limsup_{m \rightarrow \infty} \mu_{\mu^\ast h}(\widetilde{X},
                                      \mu^\ast(H+mL)) = \infty \Longrightarrow
   \limsup_{m \rightarrow \infty} \mu_h(X, H+mL) = \infty. \]
\begin{proof}
First,
$K_{\widetilde{X}} = \mu^\ast K_X + E$.
By the functorial property of multiplier ideal sheaves \cite[(5.8)]{Dem00}
this implies
$\mathcal{I}(h) = \mu_\ast (\mathcal{O}(E) \otimes \mathcal{I}(\mu^\ast h))$.
Since
$\mathcal{I}(\mu^\ast h) \subset 
 \mathcal{O}(E) \otimes \mathcal{I}(\mu^\ast h) \subset 
 \mathcal{K}_{\widetilde{X}}$
(the sheaf of total quotient rings), it follows
\[ \mu_\ast \mathcal{I}((\mu^\ast h)^{ml}) \subset \mathcal{I}(h^{ml}). \]
By the projection formula,
\[ \mu_\ast(\mu^\ast \mathcal{O}_X(l(H+mL)) \otimes 
            \mathcal{I}((\mu^\ast h)^{ml})) = \mathcal{O}_X(l(H+mL)) \otimes
   \mu_\ast \mathcal{I}((\mu^\ast h)^{ml}). \]
Consequently,
\[ h^0(\widetilde{X}, \mu^\ast \mathcal{O}_X(l(H+mL)) \otimes 
            \mathcal{I}((\mu^\ast h)^{ml}))) \leq
   h^0(X, \mathcal{O}_X(l(H+mL)) \otimes \mathcal{I}(h^{ml})), \]
which implies the claim.
\end{proof}

\noindent
\textbf{Step 3.}
$\limsup_{m \rightarrow \infty} \mu_{\mu^\ast h}(\widetilde{X},
                                      \mu^\ast(H+mL)) = \infty$.

\noindent
For
$l_0 \gg 0$
the line bundle
$l_0\mu^\ast H - E$
is ample on
$\widetilde{X}$.
It is enough to show that
\[ \limsup_{m \rightarrow \infty} \mu_{\mu^\ast h}(\widetilde{X},
                                      l_0\mu^\ast(H+mL)-E) = \infty. \]
Let
$p: \widetilde{X} \rightarrow \mathbb{P}^1$
be the projection on
$\mathbb{P}^1$.
Now, the sheaf
$\mathcal{O}_{\widetilde{X}}(ll_0\mu^\ast(H+mL)-lE) \otimes 
 \mathcal{I}((\mu^\ast h)^{mll_0})$
is torsion free. Since
$p$ 
is flat,
\[ p_\ast(\mathcal{O}_{\widetilde{X}}(ll_0\mu^\ast(H+mL)-lE) \otimes 
          \mathcal{I}((\mu^\ast h)^{mll_0})) \cong \mathcal{E}_{m,l} \cong
   \bigoplus_{i=1}^r \mathcal{O}(a_i) \]
is also torsion free, hence a locally free sheaf on 
$\mathbb{P}^1$. 
Here,
the
$a_i = a_i(m,l)$
and
$r = r(m,l)$
depend on
$m,l$.

\noindent
By upper semicontinuity and the Ohsawa-Takegoshi Extension Theorem,
for a general fibre 
$F$
\begin{eqnarray*} 
r(m,l) & =    & h^0(F, \mathcal{O}_F(ll_0\mu^\ast(H+mL)-lE) \otimes 
                    \mathcal{I}((\mu^\ast h)^{mll_0})_{|F}) \\
       & \geq & h^0(F, \mathcal{O}_F(ll_0\mu^\ast(H+mL)-lE) \otimes 
                       \mathcal{I}((\mu^\ast h)^{mll_0}_{|F})). 
\end{eqnarray*}
Since
$(l_0\mu^\ast(H+mL)-E)_{|F}$
is ample and
$(L^{l_0}_{|F}, h^{l_0}_{|F})$
is not numerically trivial by step 2, the induction hypothesis on
$F$
implies
\begin{eqnarray} \label{Edeg-eq}
\limsup_{m \rightarrow \infty}(\limsup_{l \rightarrow \infty} 
                                            l^{-(n-1)}r(m,l)) = \infty.
\end{eqnarray}

\noindent
Let
$h_0$
be a 
$\mathcal{C}^\infty$
hermitian metric on the ample line bundle
$\mathcal{O}_{\widetilde{X}}(l_0\mu^\ast H - E)$
with
$\Theta_{h_0} > 0$,
let
$h_1$
be any 
$\mathcal{C}^\infty$
metric on
$\mathcal{O}_{\mathbb{P}^1}(1)$.
Then there exists a
$c \in \mathbb{Q}_{>0}$
such that
$\Theta_{h_0} - c \pi^\ast \Theta_{h_1}$
is a positive K\"ahler form on 
$\widetilde{X}$.

\bigskip 

\noindent
\textit{Claim.}
$\mathcal{E}_{m,l} \otimes \mathcal{O}_{\mathbb{P}^1}(-cl+1)$
is globally generated for all
$l \in \mathbb{N}$
with
$lc \in \mathbb{N}$,
$l \gg 0$.
\begin{proof}
By looking at the short exact sequence 
\[ 0 \rightarrow \mathcal{E}_{m,l} \otimes \mathcal{O}_{\mathbb{P}^1}(-cl)
     \rightarrow \mathcal{E}_{m,l} \otimes \mathcal{O}_{\mathbb{P}^1}(-cl+1) 
     \rightarrow \mathcal{E}_{m,l} \otimes \mathcal{O}_{\mathbb{P}^1}/
                                                       \hbox{\gothic{m}}_x
     \rightarrow 0, \]
one sees that the vector bundle
$\mathcal{E}_{m,l} \otimes \mathcal{O}_{\mathbb{P}^1}(-cl+1)$
is globally generated if
$H^1(\mathbb{P}^1, \mathcal{E}_{m,l} \otimes \mathcal{O}_{\mathbb{P}^1}(-cl))
 = 0$.
But this cohomology group is contained in
\[ H^0(\mathbb{P}^1, 
   R^1 p_\ast(\mathcal{O}_{\widetilde{X}}(ll_0\mu^\ast(H+mL)-lE) \otimes 
              \mathcal{I}((\mu^\ast h)^{mll_0})) \otimes 
              p^\ast  \mathcal{O}_{\mathbb{P}^1}(-cl)). \]
This higher direct image sheaf is 
$0$
by Nadel vanishing \cite[(5.11)]{Dem00}, applied on preimages in
$\widetilde{X}$
of open affine subsets of
$\mathbb{P}^1$ 
and the big line bundle
$ll_0\mu^\ast(H+mL)-lE + p^\ast  \mathcal{O}_{\mathbb{P}^1}(-cl)$
equipped with the positive singular hermitian metric 
$h_0^l \otimes (p^\ast h_1)^{cl} \otimes (\mu^\ast h)^{mll_0}$.
\end{proof}

\noindent
The claim implies 
$\limsup_{l\rightarrow \infty} l^{-1} (\min_i a_i) \geq c$,
hence
\begin{eqnarray*} 
\lefteqn{l^{-n} \cdot h^0(\widetilde{X},
                 \mathcal{O}_{\widetilde{X}}(ll_0\mu^\ast(H+mL)-lE) \otimes 
                 \mathcal{I}((\mu^\ast h)^{mll_0})) = } \\
 & = & l^{-n} \cdot h^0(\mathbb{P}^1, \mathcal{E}_{m,l})\ \geq\ 
     l^{-(n-1)} r(m,l) l^{-1} (\min_i a_i)\ \geq\ c \cdot l^{-(n-1)} r(m,l) 
\end{eqnarray*}
(\ref{Edeg-eq}) implies step 3, and Lemma~\ref{vol-lem} is proven.

\begin{lem} \label{sect-lem}
Let
$X$
be a complex projective variety, let
$(L,h)$
be a pseudo-effective line bundle with positive singular hermitian metric 
$h$.
Assume that 
$(L,h)$
is not numerically trivial. Let
$x \in X$
be a sufficiently general point such that
$\mathcal{I}(h^m)_x \cong \mathcal{O}_{X,x}$.
Then, for any ample line bundle
$H$,
for all
$N \in \mathbb{N}$ 
there exists 
$m_0 \in \mathbb{N}$
such that for 
$l \gg 0$
arbitrarily large there is a section
\[ 0 \not\equiv \sigma_l \in H^0(X, \mathcal{O}_X(l(H+m_0L)) \otimes 
                \mathcal{I}(h^{m_0l}) \otimes \hbox{\gothic{m}}^{Nl}_{X,x}). \]
\end{lem}
\begin{proof}
By Lemma~\ref{vol-lem} there exists an
$m_0$
such that
the volume
$\mu_h(X, H + m_0L) > N^{\dim X} + 1$.
Consequently, for
$l \gg 0$ 
arbitrarily large,
\[ h^0(X, \mathcal{O}_X(l(H+m_0L)) \otimes \mathcal{I}(h^{m_0l})) \geq
   \frac{N^{\dim X}+1}{(\dim X)!} l^{\dim X} + o(l^{\dim X}). \]
Set
$\mathcal{F} := \mathcal{O}_X(l(H+m_0L)) \otimes \mathcal{I}(h^{m_0l})$.
Since
$\mathcal{I}(h^m)_x \cong \mathcal{O}_{X,x}$,
it is true that
$h^0(X, \mathcal{F} \otimes \mathcal{O}_X/\hbox{\gothic{m}}^{Nl}_x) =
 \frac{N^{\dim X}}{(\dim X)!} l^{\dim X} + o(l^{\dim X})$.
Using the sequence
\[ 0 \rightarrow H^0(X, \mathcal{F} \otimes \hbox{\gothic{m}}^{Nl}_x)
     \rightarrow H^0(X, \mathcal{F}) \rightarrow
     H^0(X, \mathcal{F} \otimes \mathcal{O}_X/\hbox{\gothic{m}}^{Nl}_x) \]
one gets the lemma.
\end{proof}

\subsection{The Key Lemma} \label{keylem-subsec}

\noindent
The proof of the Key Lemma~\ref{Key-lem} starts with the blow up 
$\pi: \widehat{M} \rightarrow M$
in
$W$.
Then a very general curve
$R$
in the smooth exceptional divisor
$\widehat{W} \subset \widehat{M}$
is 
$(\pi^\ast L, \pi^\ast h)-\!$
general: If
$D$
is a prime divisor in some
$E_c(h)$
with
$D \cap W \neq \emptyset$,
then the strict transform 
$\widehat{D}$
of
$D$
will have minimal Lelong number
$\nu(\pi^\ast h, \widehat{D}) \geq \nu(h, D)$.
Now choose a very general curve
$R \subset \widehat{W}$
such that the branching locus
$\pi_{|R}: R \rightarrow W$
does not contain any of the countably many points
$y \in W$
with
$\nu(h_{|W},y) > 0$.
Then for
$x \in \widehat{D}$, 
\[ \nu(\pi^\ast h_{|R},x) = \nu(h_{|W},\pi(x)) = \nu(h,y) = \nu(h,D) \leq
   \nu(\pi^\ast h, \widehat{D} \leq \nu(\pi^\ast h,x), \]
hence
$\nu(\pi^\ast h_{|R},x) = \nu(\pi^\ast h,x)$.
For all the other
$x \in R$,
the Lelong number
$\nu(\pi^\ast h_{|R},x) = 0$,
hence
\[ 0 = \nu(\pi^\ast h,x) \leq \nu(\pi^\ast h_{|R},x) = 0. \]

\noindent
Now assume that
$(L,h)$
is not numerically trivial on
$M$.
By birational invariance,
$(\pi^\ast L, \pi^\ast h)$
is not numerically trivial on
$\widehat{M}$.
(This is an application of the definition of 
$(L,h)-\!$
intersection numbers). For
an ample line bundle 
$H$
on
$\widehat{M}$, 
it follows
\[ \limsup_{m \rightarrow \infty} \mu_h(\widehat{M}, H + m\pi^\ast L) = 
   \infty \]
by Lemma~\ref{vol-lem}. Let
$x_0 \in \widehat{W}$
be a sufficiently general point such that 
$\mathcal{I}(\pi^\ast h^m)_{x_0} \cong \mathcal{O}_{\widehat{M},x_0}$
for all integers
$m$.
By Lemma~\ref{sect-lem}, for all
$N$
there exists an
$m_0$
such that for arbitrarily large
$l \gg 0$
there is a non-vanishing section
\[ \sigma_l \in H^0(\widehat{M}, \mathcal{O}_{\widehat{M}}(l(H+m_0\pi^\ast L)) 
                \otimes 
                \mathcal{I}(\pi^\ast h^{m_0l}) \otimes 
                \hbox{\gothic{m}}^{Nl}_{x_0}) - \{0\}. \]

\noindent
Let
$\mathcal{R}$
be a family of smooth intersection curves of
$n-2$
divisors in
$|H_{|\widehat{W}}|$
through
$x_0$
which cover
$\widehat{W}$. 
Choose
$d_0 \gg 0$
such that for general fibres 
$F$
of
$\widehat{f} = f \circ \pi$,
\[ H^{\dim F - 1} \cdot F \cdot (H - d_0 \widehat{W}) < 0. \]

\noindent
\textit{Claim.}
There exists an
$A_0 > 0$
independent of
$m_0$
such that 
\begin{eqnarray} \label{uplim-eq}
 & \dim H^0(R, \mathcal{O}_R(l(H + m_0 \pi^\ast L) - s\widehat{W}) \otimes 
            \mathcal{I}(\pi^\ast h^{m_0l})_{|R})
 \leq A_0 \cdot l + o(l) &
\end{eqnarray}
for all curves
$R \in \mathcal{R}$
with
$h_{|R} \not\equiv \infty$
and for all 
$0 \leq s \leq d_0l$.
\begin{proof}
Since 
$(\pi^\ast L, \pi^\ast h)$
is numerically trivial on
$\widehat{W}$,
\[ (\pi^\ast L, \pi^\ast h).R = 
   \limsup_{m \rightarrow \infty} \frac{1}{m} \deg_R 
   (\mathcal{O}_R(m\pi^\ast L) \otimes \mathcal{I}(\pi^\ast h^m)_{|R}) = 0. \]
(This is the application of Theorem~\ref{Inteq-thm}, i.e.\/ the
$(L,h)-\!$
generality of
$R$.) 
Consequently,
\begin{eqnarray*}
\deg_R (\mathcal{O}_R(l(H+m_0 \pi^\ast L) - s\widehat{W}) \otimes 
        \mathcal{I}(\pi^\ast h^{m_0l})_{|R}) & = &
\deg_R (\mathcal{O}_R(lH - s\widehat{W}) + \\
 & &  \deg_R (\mathcal{O}_R(lm_0 \pi^\ast L) \otimes 
              \mathcal{I}(\pi^\ast h^{m_0l})_{|R}) \\
 & \leq & A_0 \cdot l + o(l)\ \ 
\end{eqnarray*}
for some
$A_0 > 0$,
as long as 
$0 \leq s \leq d_0l$.

\noindent
If
$\widehat{W}.R \leq 0$ 
the ampleness of
$H$
will imply by Riemann-Roch that
\[ H^1(R, \mathcal{O}_R(l(H + m_0\pi^\ast L) - s\widehat{W}) \otimes 
          \mathcal{I}(\pi^\ast h^{m_0l})_{|R}) 
   = 0 \]
and hence the claim.

\noindent
If
$\widehat{W}.R > 0$
there will exist an
$a_0$
such that for all
$s \geq a_0$
the cohomology group
$H^0(R,  \mathcal{O}_R(s\widehat{W})) \neq 0$.
Therefore,
\begin{eqnarray*}
\lefteqn{h^0(R, \mathcal{O}_R(l(H + m_0\pi^\ast L) - s\widehat{W}) \otimes 
          \mathcal{I}(\pi^\ast h^{m_0l})_{|R}) \leq} \\ 
  & & h^0(R, \mathcal{O}_R(l(H + m_0\pi^\ast L)) \otimes 
               \mathcal{I}(\pi^\ast h^{m_0l})_{|R}),
\end{eqnarray*}
and this gives the claim for
$a_0 \leq s \leq d_0l$.
For
$s \leq a_0$
one can argue as above with 
$H$
ample and Riemann-Roch.
\end{proof}

\noindent
Now choose
$N > A_0 + d_0$.
Then 
\[ \deg_R \sigma_{l|R} \geq N \cdot l \]
for the corresponding section
\[ \sigma_{l|R} \in H^0(R, \mathcal{O}_R(l(H+m_0\pi^\ast L)) \otimes 
          \mathcal{I}(\pi^\ast h^{m_0l})_{|R} \otimes 
          \hbox{\gothic{m}}^{Nl}_{R,x_0}). \]
Because   
$s = 0$
(\ref{uplim-eq}) implies that
$\sigma_{l|R} \equiv 0$
for
$l \gg 0$
depending on
$R$. 
But since 
$\sigma_l$
vanishes on a Zariski closed subset and the curves in
$\mathcal{R}$
cover
$\widehat{W}$, 
there exists an
$l \gg 0$
such that
$\sigma_{l|\widehat{W}} \equiv 0$
and
\[ \sigma_l \in H^0(\widehat{M}, \mathcal{O}_{\widehat{M}}(l(H+m_0 \pi^\ast L)-
         \widehat{W}) 
         \otimes 
         \mathcal{I}(\pi^\ast h^{m_0l}) \otimes 
         \hbox{\gothic{m}}^{Nl-1}_{x_0}). \]
By repeating the argument for
$0 < s \leq d_0l$
one finally gets
\[ \sigma_l \in H^0(\widehat{M}, \mathcal{O}_{\widehat{M}}(l(H+m_0 \pi^\ast L)-
         d_0l\widehat{W}) 
         \otimes \mathcal{I}(\pi^\ast h^{m_0l})).\]

\noindent
Let
$F$
be sufficiently general,
$\pi^\ast h_{|F} \not\equiv \infty$
and
$(\pi^\ast L,\pi^\ast h)$
numerically trivial on
$F$.
Let
$\mathcal{S}_F$
be a family of smooth intersection curves of
$\dim F - 1$
divisors in
$|H_{|F}|$
covering
$F$. 
Let
$S \in \mathcal{S}_F$
be such a curve, 
$\pi^\ast h_{|S} \not\equiv \infty$,
and assume that
$S \cap \widehat{W} \neq \emptyset$.
Since 
$(\pi^\ast L, \pi^\ast h)$
is numerically trivial on
$F$,
\[ (\pi^\ast L, \pi^\ast h).S = \limsup_{m \rightarrow \infty} \frac{1}{m} 
   \deg_S 
   (\mathcal{O}_S(m\pi^\ast L) \otimes \mathcal{I}(\pi^\ast h^m)_{|S}) = 0. \]
Furthermore, by assumption
\[ S.(H-d_0\widehat{W}) = H^{\dim F - 1}.F.(H-d_0\widehat{W}) < 0, \]
hence
\[ \deg_S 
   (\mathcal{O}_S(l(H+m\pi^\ast L)-d_0l\widehat{W}) \otimes 
    \mathcal{I}(\pi^\ast h^{m_0l})_{|S}) < 0 \]
for some 
$l \gg 0$,
and as above one concludes
$\sigma_{l|F} \equiv 0$,
$\sigma_l \equiv 0$
which is a contradiction.

\subsection{Proof of the pseudo-effective Reduction Map Theorem}

The main construction used in this proof is described by the following
\begin{lem} \label{famconst-lem}
Let
$X$
be a complex projective variety, let
$M$
be a set of subvarieties 
$F_m \subset X$,
$m \in M$,
such that the union
$\bigcup_{m \in M} F_m \subset X$
is not contained in a pluripolar set in 
$X$.
Then there is a family
$\mathfrak{F} \subset X \times B$
of subschemes of
$X$,
covering the whole of
$X$,
and a set
$B^\prime \subset B$
not contained in a pluripolar set of
$B$,
parametrizing subvarieties
$F_m$,
$m \in M$.
\end{lem}
\begin{proof}
$M$ 
may be interpreted as a subset of 
$\mathrm{Chow(X)}$. 
There are only countably many components of
$\mathrm{Chow(X)}$.
Hence there must be at least one component
$\mathcal{C} \subset \mathrm{Chow(X)}$
such that the subschemes parametrized by the Zariski closure 
$Z = \overline{\mathcal{C} \cap M}$
cover the whole of
$X$,
and 
$\mathcal{C} \cap M$
is not a pluripolar set in
$Z$.
Otherwise, the subvarieties
$F_m$,
$m \in M$,
are contained in a pluripolar set of
$X$,
contradiction.
\end{proof}

\noindent
Consider families
$\widetilde{f}: \mathfrak{X} \rightarrow \mathcal{N}$
with the following properties:
\begin{itemize}
\item[(i)]
$\mathfrak{X} \subset X \times \mathcal{N}$,
$\mathfrak{X}, \mathcal{N}$
quasi-projective, irreducible, general fibres are subvarieties of
$X$;

\smallskip

\item[(ii)]
the projection
$p: \mathfrak{X} \rightarrow X$
is generically finite;

\smallskip

\item[(iii)]
$(L,h)$
is defined and numerically trivial on sufficiently general fibres of 
$\widetilde{f}$,
i.e. on a set of fibres
$\mathcal{M} \subset \mathcal{N}$
which is not contained in a pluripolar set;

\smallskip

\item[(iv)]
the fibres are generically unique, i.e. if
$U \subset \mathcal{N}$
is an open subset such that
$\widetilde{f}_{|U}$
is flat then the induced map 
$U \rightarrow \mathrm{Hilb}(X)$
will be generically bijective.
\end{itemize}
The identity map
$\mathrm{id}: X \rightarrow X$
is such a family, hence there is one with minimal base dimension 
$\dim \mathcal{N}$.

\smallskip

\noindent
\textit{Claim.}
The projection
$p: \mathfrak{X} \rightarrow X$
is birational on such a minimal family
$\widetilde{f}: \mathfrak{X} \rightarrow \mathcal{N}$.
\begin{proof}
Assume that
$p$
is not birational.

\noindent
Then, for a general fibre
$F$
of
$f$
and a general point
$x \in F$
there is another fibre
$F^\prime$
containing 
$x$,
hence a curve
$C^\prime$
with
$x \in C^\prime \subset F^\prime$
and
$C^\prime \not\subset F$.
Consequently, one gets a family of curves
$g: \mathfrak{C} \rightarrow \mathcal{M}$
with
$\mathfrak{C} \subset \mathfrak{X} \times \mathcal{M}$
giving a generically finite covering of
$\mathfrak{X}$
such that for the general 
$g-\!$
fibre curve
$C$
the 
$f-\!$
projection has
$\dim f(C) = 1$.
By blowing up and base change one can assume the following situation:
\[ \xymatrix{ 
   X & {\mathfrak{X}} \ar[l]_{p} \ar[d]^{f} & 
       {\widetilde{\mathfrak{C}}} \ar[l]_{\pi} \ar[d]^{\widetilde{f}} 
                                  \ar[r]^{\widetilde{g}}&
   {\mathcal{\widetilde{M}}}  \\
     & {\mathcal{N}} & {\mathcal{\widetilde{N}}} \ar[l] &
             } \]
where
$\widetilde{\mathfrak{C}}$,
$\mathcal{\widetilde{M}}$
and
$\mathcal{\widetilde{N}}$
are smooth, the general fibres of
$\widetilde{g}$
are smooth curves and the fibres of
$\widetilde{f}$
map onto fibres of
$f$
in 
$X$.
Furthermore, the maps
$p$
and
$\pi$
are generically finite.

\noindent
Let
$(\widetilde{L},\widetilde{h})$
be the pulled back 
$(L,h)$. 
Take an
$(\widetilde{L},\widetilde{h})-\!$
general
$\widetilde{g}-\!$
fibre curve
$C$
in
$\widetilde{\mathfrak{C}}$
such that the general
$\widetilde{f}-\!$
fibre through points of
$C$
is smooth. Look at the subvariety
$G_C = \widetilde{f}^{-1}(\widetilde{f}(C)) \subset \widetilde{\mathfrak{C}}$.
It may be not smooth, but by the smoothness of the general
$\widetilde{f}-\!$
fibre, the singular locus does not contain 
$C$.
Hence using Lemma~\ref{birinvLh-lem}, an embedded resolution of
$G_C$
in
$\widetilde{\mathfrak{C}}$
gives a smooth subvariety
$\widehat{G_C}$
in the blow up
$\widehat{\mathfrak{C}}$
such that the strict transform 
$\widehat{C}$
of
$C$
is still 
$(\mu^\ast \widetilde{L}, \mu^\ast \widetilde{h})-\!$
general in
$\widehat{\mathfrak{C}}$.
By the following lemma,
$\widehat{C}$
is also 
$(\mu^\ast \widetilde{L}, \mu^\ast \widetilde{h})-\!$
general in
$\widehat{G_C}$, 
and one can apply the Key Lemma:
$\widehat{G_C}$
is 
$(\mu^\ast \widetilde{L}, \mu^\ast \widetilde{h})-\!$
numerically trivial. By birational invariance this is true for the image of
$\widehat{G_C}$
in
$X$,
too. But 
$\dim \widehat{G_C} > \dim F$.
Since all curves in a family are 
$(L,h)-\!$
general outside a pluripolar set,
the construction in Lemma~\ref{famconst-lem} 
gives a new family
\[ g: \mathfrak{Y} \rightarrow \mathcal{M} \]
satisfying conditions (i) - (iv), and 
\[ \dim \mathcal{M} = \dim X - \dim G = \dim X - (\dim \mathcal{N}^\prime +
   \dim F) < \dim X - \dim F = \dim \mathcal{N} \]
for fibres
$G$
of
$g$.
This is a contradiction to the minimality of
$\dim \mathcal{N}$.
\end{proof}

\begin{lem}
Let
$Y \subset X$
be a smooth subvariety in a projective complex variety
$X$
with a pseudo-effective line bundle
$L$
and a positive singular hermitian metric
$h$
on
$L$
such that
$h_{|Y} \not\equiv \infty$.
Then an
$(L,h)-\!$
general curve is also an
$(L_{|Y},h_{|Y})-\!$
general curve on 
$Y$.
\end{lem}
\begin{proof}
$\nu(h_{|C},x) = 0$
implies
$0 =\nu(h_{|Y},x) \leq \nu(h_{|C},x) = 0$,
hence
$\nu(h_{|Y},x) = \nu(h_{|C},x)$.

\noindent
$\nu(h_{|C},x) > 0$
implies that
$x \in D$
for some prime divisor
$D$
on some
$E_c(h)$.
The restricted divisor 
$D_{|Y}$
may be singular but not in 
$x$:
Then,
\[ \nu(h_{|Y},x) > \nu(h_{|Y},D_{|Y}) \geq \nu(h,D) = \nu(h_{|C},x), \]
contradiction.
\end{proof}

\noindent
In the same way one shows

\noindent
\textit{Claim.} Let
$\widetilde{g}: \mathfrak{X} \rightarrow \widetilde{\mathcal{N}}$
be another family satisfying the conditions (i) - (iv). Then there is a 
commutative diagram of rational maps
\[ \xymatrix{
   X & {\mathfrak{X}} \ar[l] \ar[d]^{\widetilde{f}} & 
       {\widetilde{\mathfrak{X}}} \ar@{-->}[l] \ar[d]^{\widetilde{g}} \\
     & {\mathcal{N}} & {\widetilde{\mathcal{N}}} \ar@{-->}[l]
            } \]
such that the general fibre of
$\widetilde{g}$
is contained in a fibre of
$\widetilde{f}$.
\hfill $\Box$

\noindent
On the one hand, this claim implies the birational uniqueness of 
$\widetilde{f}$. 
On the other hand one can prove claim (ii) in the pseudo-effective Reduction 
Map Theorem: If (ii) is not satisfied there will be a set of points
$N \subset X$
which is not contained in a pluripolar set such that
\[ \forall\ x \in N\ \exists\ C_x \ni x\ \mathrm{irreducible\ curve,\ } \dim 
   \widetilde{f}(C_x) = 1: (L,h).C_x = 0. \]
By Lemma~\ref{famconst-lem} one gets a family of curves satisfying conditions
(i)-(iv). The claim implies that the general fibre of this family 
is contained in a fibre of
$\widetilde{f}$,
hence also some of the curves
$C_x$:
contradiction.

\noindent
Finally it is possible to prove that in part (i) of the Reduction Map Theorem,
all fibres outside a pluripolar set are 
$(L,h)-\!$
numerically trivial: This pluripolar set is just the set of fibres 
$F$
such that
$h_{|F} \equiv \infty$.
Because assume to the contrary that
$C \subset F$
is a curve on a fibre 
$F$
such that
$h_{|C} \not\equiv \infty$,
hence
$h_{|F} \not\equiv \infty$,
and 
$C$
is not
$(L,h)-\!$
numerically trivial. Then, as in step 1 of subsection~\ref{vol-ssec},
$(L,h)$
is not numerically trivial on sufficiently general fibres
$F$,
contradiction !

\newcommand{\etalchar}[1]{$^{#1}$}

\end{document}